\documentclass{amsart}
\usepackage{tabu}
\usepackage{amssymb} 
\usepackage{amsmath} 
\usepackage{amscd}
\usepackage{amsbsy}
\usepackage{comment, enumerate}
\usepackage[matrix,arrow]{xy}
\usepackage{hyperref}
\usepackage{mathrsfs}

\DeclareMathOperator{\ord}{ord}

\newcommand{\Q}{{\mathbb Q}}
\newcommand{\Z}{{\mathbb Z}}

\newcommand{\OO}{{\mathcal O}}
\newcommand{\ga}{{\mathfrak{a}}}

\newcommand{\fp}{\mathfrak{p}}
\newcommand{\fq}{\mathfrak{q}}

\newcommand{\fz}{\mathfrak{z}}

\begin {document}

\newtheorem{thm}{Theorem}
\newtheorem{lem}{Lemma}[section]

\theoremstyle{definition}

\theoremstyle{remark}

\title[]{Perfect powers that are sums of consecutive squares}


\author{Vandita Patel}
\address{Department of Mathematics, University of Toronto, Bahen Centre, 40 St. George St., Room 6290, Toronto, Ontario, Canada, M5S 2E4}
\email{vandita\_patel@hotmail.co.uk}

\date{\today}


\begin {abstract}
We determine all perfect powers that can be written as the sum of at most
$10$ consecutive squares.
\end {abstract}
\maketitle

\section{Introduction} \label{intro}

Power values of power sums of consecutive integers have been of interest throughout the past 70 years. Techniques from algebraic number 
theory and Diophantine approximation have allowed the resolution of such equations with small exponents, as well as proofs of general
theorems. 
This can be seen in the work of Brindza \cite{Br}, Cassels \cite{Cassels}, Gy\H{o}ry, Tijdeman and  Voorhoeve \cite{GTV}, 
Hajdu \cite{Ha}, Pint\'er \cite{P}, Sch\"affer \cite{Sc}, and Zhang and Bai \cite{ZB} among many others.

Recently, Bennett, Gy\H{o}ry and Pint\'er \cite{BGP}, Pint\'er \cite{P2}, Zhang \cite{Zhang},
 as well as the author in collaboration with Bennett and Siksek \cite{BPS1}, \cite{BPS2} resolved many such equations using the modular method. A Galois-theoretic approach
can be found in \cite{PatelSiksek}. 

\medskip

In this paper we consider the equation
\begin{equation} \label{eqn:main}
(x+1)^2+(x+2)^2+\cdots+(x+d)^2 =y^n, \qquad n \ge 2,
\end{equation}
with $2 \le d \le 10$. Here $x$ and $y$ denote integers. This equation was solved by Cohn \cite{Cohn}
for $d=2$ and by Zhang \cite{Zhang} for $d=3$.
\begin{thm}\label{thm:main}
Let $2\leq d \leq 10$. The only solutions to equation \eqref{eqn:main} with $n \ge 3$
are 
\begin{gather*}
(d,x,y,n)=(2,-1,\pm 1,2r), \quad (2,-2,\pm 1,2r), \quad 
(2,-1, 1,2r+1), \\
(2,-2,1,2r+1), \quad 
(2, 118 , \pm 13 ,4), \quad (2,-121,\pm 13, 4).
\end{gather*}
The only solutions to equation \eqref{eqn:main} with $n=2$ is the infinite family with $d=2$ and
\[
2x+3+y \sqrt{2}=\pm (1+\sqrt{2})^{2r+1}, \qquad r \in \Z.
\]
\end{thm}

\section{Sums of Consecutive Squares}

We make use of the following Lemma of Zhang and Bai \cite{ZB}.
\begin{lem}[Zhang and Bai]\label{lem:ZB}
Let $p$ be a prime such that $ p \equiv \pm5 \mod 12$. If $p \mid d$ and $\ord_p(d) \not\equiv 0 \mod n$ then equation~\eqref{eqn:main} has no solutions.
\end{lem}
We note that Lemma~\ref{lem:ZB} immediately allows us to eliminate the cases $d=5$, $7$ and  $10$. As mentioned in the introduction,
equation~\eqref{eqn:main} has been solved by Zhang \cite{Zhang} for $d=3$. It was solved for $d=2$ and $n \ge 3$ by 
Cohn \cite{Cohn}, and this gives the solutions enumerated in the theorem for $n \ge 3$. 

It remains to consider $d=2$ with $n=2$, and $d=4$, $6$, $8$, $9$ with $n \ge 2$.
We rewrite equation~\eqref{eqn:main}  as
\begin{equation}\label{eqn:conssquares1}
dx^2 + d(d+1)x + \frac{d(d+1)(2d+1)}{6} = y^n.
\end{equation}

Factorising and completing the square gives us
\begin{equation}\label{eqn:conssquares}
d \left(\left( x + \frac{d+1}{2} \right)^2  + \frac{(d-1)(d+1)}{12} \right) = y^n.
\end{equation}

\begin{lem}\label{lem:1}
Let $r = \ord_2(d)$. Suppose $r\geq 2$, then in equation~\eqref{eqn:conssquares} we have $n \mid (r-1)$.
\end{lem}
\begin{proof}
Let $D = d/2^2$. We substitute into equation~\eqref{eqn:conssquares} to get,
\begin{equation}\label{eqn:dval}
D \left(\left( 2x + (d+1) \right)^2  + \frac{(d-1)(d+1)}{3} \right) = y^n.
\end{equation}
Observe that 
\[
(2x+d+1)^2 \equiv 1 \pmod{4}, \qquad
\frac{(d-1)(d+1)}{3} \equiv 1 \pmod{4}.
\]
Comparing valuations on both sides of \eqref{eqn:dval} we see that
\[
n \ord_2(y)=\ord_2(D)+1=r-1.
\]
This completes the proof.
\end{proof}

\begin{lem}\label{lem:2}
Let $q = 3$ and  $r = \ord_3(d)$. Suppose $r\geq 2$,  then in equation~\eqref{eqn:conssquares}, we have $n \mid (r-1)$.
\end{lem}
\begin{proof}
Let $D = d/3$. We substitute into equation~\eqref{eqn:conssquares} to get,
\[
D \left(3\left( x + \frac{(d+1)}{2} \right)^2  + \frac{(d-1)(d+1)}{4} \right) = y^n.
\]
Observe that the expression in brackets is never divisible by 3. Hence 
$\ord_3(D)=\ord_3(y^n)=n \ord_3(y)$, proving that
$n \mid (r-1)$.
\end{proof}
Applying Lemmata~\ref{lem:1} and~\ref{lem:2} allows us to eliminate $d=4$ and $d=9$,
and $d=8$  with $n \ge 3$. For the proof of Theorem~\ref{thm:main}, it remains to deal
with $d=6$, and also with $d=2$, $8$ for $n=2$.

\section{Case: $n=2$}

In this section, we deal with equation~\eqref{eqn:main} with $n=2$ and $d=2$, $6$, $8$.
First, we consider $d=2$. Then, equation \eqref{eqn:main} can be rewritten as
\[
(2x+3)^2-2y^2=-1.
\]
This yields the infinite family of solutions in Theorem~\ref{thm:main}.

Now let $d=6$; we can rewrite equation~\eqref{eqn:main} as 
\[
3(2x+7)^2+35=2y^2
\]
which is impossible as the left-hand side is $6 \pmod{8}$ and the right-hand side is $2 \pmod{8}$.

Finally let $d=8$; we can rewrite equation~\eqref{eqn:main} as
\[
2\left((2x+9)^2+21 \right)=y^2.
\]
Write $y=2Y$, we obtain
\[
(2x+9)^2+21=2 Y^2.
\]
Again, we see that the left-hand side is $6 \pmod{8}$, yielding a contradiction.

\section{The Case $d=6$}

It finally remains to solve equation~\eqref{eqn:main} for $d=6$. We suppose $n=p$ is an odd prime.
We rewrite equation~\eqref{eqn:main} as 
\begin{equation}\label{eq:6squares}
X^2 + 3\times 5\times 7 = 6y^p,
\end{equation}
where $X=6x+21$.  It is easy to see that $2,3,5,7 \nmid y$.
Let $K = \Q(\sqrt{-105})$ and write $\OO_K=\Z[\sqrt{-105}]$ for its ring of integers.
This has class group isomorphic to $(\Z/2\Z)^3$. We factorise the left-hand side
of equation~\eqref{eq:6squares} as
\[
(X + \sqrt{-105})(X - \sqrt{-105}) = 6y^p.
\]
It follows that
\[
(X + \sqrt{-105})\OO_K = \fp_2\fp_3 \cdot \fz^p
\]
where $\fp_2$ and $\fp_3$ are the unique primes of $\OO_K$ above $2$ and $3$ respectively, and $\fz$
is an ideal of $\OO_K$.
Let $\ga = \fp_2\fp_3$. Then
\[
(X + \sqrt{-105})\OO_K = \fp_2\fp_3 \cdot \fz^p =\ga^{1-p} \cdot (\ga \fz)^p 
=(6^{(1-p)/2})(\ga \fz)^p.
\]
It follows that $\ga \fz$ is a principal ideal. Write $\ga\fz=\gamma \OO_K$ where $\gamma \in \OO_K$,
and $\ord_{\fp_2}(\gamma)=\ord_{\fp_3}(\gamma)=1$.
After possibly changing the sign of $\gamma$ we obtain,
\begin{equation}\label{eqn:X}
X+\sqrt{-105}=\frac{\gamma^p}{6^{(p-1)/2}}.
\end{equation}
Subtracting the conjugate equation from this equation, we obtain 
\begin{equation}\label{eq:cons6pp2}
\frac{\gamma^p}{6^{(p-1)/2}} - \frac{\bar{\gamma}^p}{6^{(p-1)/2}} = 2\sqrt{-105}.
\end{equation}

Let $L = \Q(\sqrt{-105}, \sqrt{6}) =  \Q(\sqrt{-70}, \sqrt{6}) $. 
Write $\OO_L$ for the ring of integers of $L$ and let
\[
\alpha = \frac{\gamma}{\sqrt{6}} \qquad \text{and} \qquad 
\beta = \frac{\bar{\gamma}}{\sqrt{6}}.
\]
Substituting into equation~\eqref{eq:cons6pp2}, we see that
\begin{equation}\label{eqn:70}
\alpha^p - \beta^p = \sqrt{-70}.
\end{equation}
\begin{lem}\label{lem:lehmer}
Let $\alpha$, $\beta$ be as above. Then
$\alpha$ and $\beta$ are algebraic integers. Moreover, $(\alpha + \beta)^2$ and $\alpha\beta$ are non-zero, coprime, rational integers,
and $\alpha/\beta$ is not a root of unity.
\end{lem}
\begin{proof}
Note that $\ga \cdot \OO_L=\sqrt{6} \OO_L$. As $\ga \mid \gamma$, $\overline{\gamma}$, we have
$\alpha$, $\beta \in \OO_L$.

Now $\gamma=u+v \sqrt{-105}$ with $u$, $v \in \Z$. Thus 
\[
(\alpha+\beta)^2=\frac{2u^2}{3}.
\]
but $\fp_3 \mid \gamma$ and $\fp_3 \mid \sqrt{-105}$, hence $\fp_3 \mid u$ and so $3 \mid u$. Therefore $(\alpha+\beta)^2 \in \Z$.
If $(\alpha+\beta)^2=0$ then $u=0$ and from equation~\eqref{eqn:X}, we establish that $X=0$, which doesn't result in an integer solution since $X=6x+21$. 
Therefore, $(\alpha+\beta)^2$ is a non-zero rational integer.

Furthermore, $\alpha \beta=(\gamma \overline{\gamma})/6$, which is clearly a non-zero rational integer.

We must check that $(\alpha+\beta)^2$ and $\alpha \beta$ are coprime. Suppose they are not coprime. 
Then there is some prime ideal $\fq$ of $\OO_L$ dividing both. This divides $\alpha$, $\beta$, and so
using  equation~\eqref{eqn:70}, we see that $\ord_\fq(\sqrt{-70}) \ge p$ and arrive at a contradiction. 

Finally we need to show that $\alpha/\beta=\gamma/\overline{\gamma} \in K=\Q(\sqrt{-105})$ is not
a root of unity. But the only roots of unity in $K$ are $\pm 1$. If $\alpha/\beta=\pm 1$
then from equation~\eqref{eqn:70}, we obtain $0=\sqrt{-70}$ or $2 \alpha^p=\sqrt{-70}$, both giving
a contradiction.
\end{proof}

Continuing with the notation of the previous proof we have, 
\[
\alpha-\beta=2v \sqrt{-105}{\sqrt{6}}=v \sqrt{-70}.
\]
Therefore, equation~\eqref{eqn:70} gives $v= \pm 1$ and 
\begin{equation}\label{eqn:thue}
\frac{\alpha^p - \beta^p}{\alpha - \beta} = \pm 1.
\end{equation}
To complete the proof, we need a famous theorem due to Bilu, Hanrot and Voutier \cite{BHV}.
Attached to a pair of algebraic numbers $\alpha$ and  $\beta$ satisfying Lemma~\ref{lem:lehmer}
is a \textbf{Lehmer sequence} given by
\[
\tilde{u}_m=
\begin{cases}
(\alpha^m-\beta^m)/(\alpha-\beta) \qquad \text{$m$ odd}\\
(\alpha^m-\beta^m)/(\alpha^2-\beta^2) \qquad \text{$m$ even}.
\end{cases}
\]
A prime $q$ is said to be a \textbf{primitive divisor} for $\tilde{u}_m$ if divides $\tilde{u}_m$
but does not divides $(\alpha^2-\beta^2)^2 \cdot \tilde{u}_1 \cdot \tilde{u}_2 \cdots \tilde{u}_{m-1}$.
The pair $(\alpha,\beta)$ is said to be \textbf{$m$-defective} if $\tilde{u}_m$ does not have a primitive divisor.
Observe that \eqref{eqn:thue} implies that $(\alpha,\beta)$ is $p$-defective. 

By Theorem 1.4 of \cite{BHV}, if $m \ge 30$ then $\tilde{u}_m$ must have a primitive divisor. Thus
we know that $p<30$. To deal with primes in the range $7 \le p <30$ we need the results of Voutier \cite{Voutier}.
The only possible values of $p$ in that range for which $\tilde{u}_p$ has no primitive divisor are $p=7$, $13$,
and for these values Voutier gives the possibilies for $\alpha/\beta$. Examining his table quickly eliminates
these cases as it is incompatible with $\alpha/\beta=\gamma/\overline{\gamma} \in \Q(\sqrt{-105})$. 
This proves the proposition.

It remains to deal with $p=3$ and $p=5$. We may rewrite equation~\eqref{eqn:70} as
\[
(u+\pm \sqrt{-105})^p-(u \mp \sqrt{-105})^p=\sqrt{6}^p \cdot \sqrt{-70}.
\]
We merely check that these polynomial equations do not have roots in $\Z$.
This completes the proof.


\end{document}